\documentclass[11pt,a4paper]{amsart}
\usepackage{amsfonts}
\usepackage{amssymb}
\usepackage[latin1]{inputenc}
\usepackage{amsmath}
\usepackage{delarray}
\newtheorem{thm}{Theorem}[section]

\newtheorem{prop}[thm]{Proposition}
\theoremstyle{definition}

\numberwithin{equation}{section}

\newcommand{\cB}{\mathcal B}

\newcommand{\cL}{\mathcal L}

\newcommand{\cU}{\mathcal U}

\newcommand{\bbR}{\mathbb R}
\newcommand{\bbT}{\mathbb T}

\newcommand{\rank}{{\rm rank\ }}

\begin{document}

\title[Kolmogorov condition]{Kolmogorov condition near hyperbolic
singularities of integrable Hamiltonian systems}

\author{Nguyen Tien Zung}
\address{Institut de Mathématiques de Toulouse, Universit\'{e} Toulouse III}
\email{tienzung.nguyen@math.ups-tlse.fr}

\date{Version 3, August 2007}

\maketitle

\begin{abstract}{

In this paper we show that, if an integrable Hamiltonian system
admits a nondegenerate hyperbolic singularity then it will satisfy
the Kolmogorov condegeneracy condition near that singularity (under
a mild additional condition, which is trivial if the singularity
contains a fixed point). \\



}
\end{abstract}


\section{Introduction}

The celebrated Kolmogorov--Arnold--Moser theorem (e.g.,
\cite{Arnold,Kolmogorov,Moser}) says that, under a small
perturbation, most invariant tori of an integrable Hamiltonian
system persist. This theorem is stated under a non-degeneracy
condition, called the Kolmogorov condition, which says that the
Hessian of the integrable Hamiltonian function $H$ with respect to a
family of action variables $(I_i)$ does not vanish: $ \det (
\partial^2 H /
\partial I_i\partial I_j ) \neq 0$. There are many generalizations of this
theorem, which require a weaker non-degeneracy condition than the
Kolmogorov condtion (see, e.g., Rüssmann \cite{Russmann-KAM2001}).
However, the Kolmogorov condition is quite natural, and integrable
Hamiltonian systems which are not resonant are expected to satisfy
this condition in general. On the other hand, in practice, this
condition is not easy to verify directly, because the computation of
the above determinant often involves Abelian integrals and
transcendental functions (see, e.g., Horozov \cite{Horozov} for the
case of spherical pendulum).

In this paper, we will show that if an integrable Hamiltonian system
admits a nondegenerate singularity of hyperbolic type, then the
Hessian $ \det ( \partial^2 H / \partial I_i\partial I_j ) \neq 0$
is different from zero everywhere in (the regular part of) a
neighborhood of that hyperbolic singularity, provided that the
integrable subsystem on a corresponding center manifold satisfies
the Kolmogorov condition (if the singularity has a fixed point, i.e.
a point at which the differential of the momentum map vanishes, then
the center manifold is just a point, and this last condition is
empty). In fact, we will show the following asymptotic formula for
$\det (\partial^2 H /
\partial I_i\partial I_j )$ near a nondegenerate hyperbolic singular
fiber $N$ of corank $k$ of the system (a \emph{fiber} means a
connected component of a level set of the momentum map; the
\emph{corank} $k$ is the maximal corank of the differential of the
momentum map on $N$):
\begin{equation}
\label{eqn:AsymptoticDet} \det ( \partial^2 H / \partial I_i\partial
I_j )(z) = {g(z) \over \prod_{i=1}^k F_{n -k + i}(z) (\ln F_{n -k +
i}(z))^3},
\end{equation}
where  $F_{n-k+1},\hdots,F_{n}$ are a well-chosen set of smooth
first integrals which vanish on $N$ (more precisely, these functions
are chosen so that the local bifurcation diagram of the momentum map
near the image of $N$  is a union of $k$ transversal hypersurfaces
given by $\prod_{i=1}^k F_{n - k + i} = 0$), $z$ represents a
regular fiber (i.e., a Liouville torus) near the singularity, and
$g(z)$ is a first integral in a connected component of a
neighborhood of $N$ minus the singular fibers, such that the limit
$\lim_{z \to N} g(z)$ exists and is different from zero. Formula
(\ref{eqn:AsymptoticDet}) implies immediately that $\det (
\partial^2 H / \partial I_i\partial I_j )(z) \neq 0$ for $z$ close
enough to the hyperbolic singular fiber $N$, and moreover $\det (
\partial^2 H / \partial I_i\partial I_j )(z) \to \infty$
when $z$ tends to  $N$.

Our result may be viewed as a significant improvement of a result of
Knörrer \cite{Knorrer}, which says that the Kolmogorov condition
$\det (\partial^2 H / \partial I_i\partial I_j ) \neq 0$ is
satisfied \emph{almost everywhere} near a nondegenerate hyperbolic
singular fiber of corank 1 or 2. Here we show that it is satisfied
\emph{everywhere} near the singular fiber, and in our result there
is no restriction on the corank of the singularity.

A similar asymptotic formula for $\det (\partial^2 H /
\partial I_i\partial I_j )$ near simple focus-focus singularities
(with just 1 singular point on the singular fiber) of an integrable
Hamiltonian system with 2 degrees of freedom was obtained by Rink
\cite{Rink-Focus2004} and Dullin and Vu-Ngoc
\cite{DullinSan-Twist2004}. We suspect that similar results hold for
any nondegenerate singularity without elliptic components, and for
many degenerate singularities as well, although the question remains
open even for the case of a focus-focus singularity with several
singular points on the singular fiber in a integrable system with 2
degrees of freedom, to my knowledge. The reason is that asymptotic
formulas for the action functions near a generic focus-focus
singularity with more than one singular points can be quite more
complicated than the case with just one singular point or the
hyperbolic case.

Most singularities of finite-dimensional integrable Hamiltonian
systems are nondegenerate in a natural sense, and a large part of
these nondegenerate singularities are of hyperbolic type. For
example, most integrable cases of rigid body problems (see, e.g.,
Chapter 14 of \cite{BolsinovFomenko-IntegrableBook}), geodesic flows
on multi-dimensional ellipsoids, finite-dimensional subsystems of
the integrable focusing cubic non-linear Schrodinger equations or
the sine-Gordon equation, etc., admit hyperbolic singularities of
various rank and corank, and our result can be applied to them.

Remark that, if the system is analytic and if $\det (\partial^2 H /
\partial I_i\partial I_j ) \neq 0$ somewhere then it is different from zero
almost everywhere, at least in a connected component of the set of
Liouville tori. So even though our result has a local character, it
can be applied to show that the Kolmogorov condition is satisfied
almost everywhere globally. On the other hand, when an integrable
Hamiltonian system is perturbed, then Liouville tori which are too
near unstable singularities are destroyed due to phenomena like
separatrix splitting (looking at it another way, a small global
perturbation will look big in action-angle coordinates in a
neighborhood of a Liouville torus which is too close to a unstable
singularity, and so KAM theory does not apply there). So the
applicability of our result to KAM theory for Liouville tori which
are very close to hyperbolic singularities is quite limited.

The rest of this paper is organized as follows:  in Section 2 we
will recall some known facts about the structure of nondegenerate
singularities of integrable Hamiltonian systems and give a more
precise statement of our main result, in Section 3 we will recall an
asymptotic formula for the action functions near nondegenerate
hyperbolic singularities, and in Section 4 we will prove the
asymptotic formula (\ref{eqn:AsymptoticDet}).

{\bf Acknowledgements}. This work is supported by the French ANR
research projects JC05-41465 (Intégrabilité réelle et complexe en
mécanique hamiltonienne) and ANR-05-BLAN-0029-01 (GIMP). A part of
this paper was written during the author's visit to Max-Planck
Institut für Mathematik, and he would like to thank MPIM for its
hospitality an excellent working conditions. I would like to thank
Vu Ngoc San for interesting discussions, the referee for his
critical remarks which helped improve the presentation of this
paper, and the editors of the special volume in honor of Richard
Cushman, who are also the organizers of a conference in his honor in
Utrecht, for the invitation to give a talk and to submit this paper.
It is a great pleasure for me to dedicate this paper to Richard
Cushman.

\section{Hyperbolic singularities}
\label{section:hyperbolic}

In order to state our result more precisely, let us recall here some
facts and definitions (see, e.g., \cite{ZungA-L,Zung-IntegrableII,
BolsinovFomenko-IntegrableBook}). Denote by ${\bf F} =
(F_1,\hdots,F_n) : (M,\omega) \rightarrow \bbR^n$ a smooth momentum
map of an integrable Hamiltonian function $H$ on a $2n$-dimensional
symplectic manifold $M$ with the symplectic form $\omega$. We will
always assume that the map $\bf F$ is \emph{proper}. Then, according
to the classical Liouville-Mineur theorem
\cite{Mineur-AA1935,Mineur-AA1937}, each connected component $T$ of
a regular level set of the momentum map $\bf F$ is an
$n$-dimensional torus, called a Liouville torus, and in a
neighborhood $\cU(T)$ of $T$ there is a so-called action-angle
coordinate system $(I_1,q_1,\hdots,I_2,q_2)$, where $q_i$ are cyclic
coordinates (defined modulo 1), such that the symplectic form is
$\omega = \sum_{i=1}^n d I_i \wedge d q_i$, and the first integrals
$F_i$ depend only on the action variables $I_1,\hdots,I_n$.

A point $y \in M$ is called a \emph{singular point} of the system if
$\rank d{\bf F} (y) < n$. The number $n - \rank d{\bf F} (p)$ is
called the \emph{corank}. If $d{\bf F} (y) = 0$ then we say that $y$
is a \emph{fixed point}. When $y$ is a fixed point, then it makes
sense to talk about the quadratic part ${\bf F}^{(2)} =
(F_1^{(2)},\hdots,F_n^{(2)})$ of the momentum map at $y$. The
functions $F_1^{(2)},\hdots,F_n^{(2)}$ are quadratic functions on
the tangent space $T_y M$, which Poisson-commute with respect to
$\omega(y)$. The space of quadratic functions on $T_yM$ together
with the Poisson bracket is naturally isomorphic to the Lie algebra
$Sp(2n,\bbR)$ of infinitesimal symplectic linear transformations,
and $(F_1^{(2)},\hdots,F_n^{(2)})$ span an Abalian subalgebra of
this Lie algebra. $y$ is called \emph{nondegenerate} if this Abelian
subalgebra is a Cartan subalgebra of $Sp(2n,\bbR)$. More generally,
a singular point $y$ of corank $k$ is called \emph{nondegenerate} if
it becomes a fixed nondegenerate singular point after a local
reduction with respect to a local free Poisson $\bbR^{n-k}$-action
generated by $(n-k)$ components of the momentum map near $y$.

According to the linearization theorem for nondegenerate singular
points, due to Vey \cite{Vey-Separable1978} in the analytic case and
Eliasson \cite{Eliasson-NF1990} in the smooth case, near a
nondegenerate singular point $y$ of corank $k$ there is a local
smooth symplectic coordinate system $(p_1,q_1,\hdots,p_n,q_n)$, such
that the functions $f_1,\hdots,f_{k},p_{k+1},\hdots,p_{n}$ are local
first integrals of the system, where each $f_i$ is of one of the
following three types:

$ \; f_i = p_i^2 + q_i^2$  (elliptic type)

$ \; f_i = p_iq_i$  (hyperbolic type)

$ \left.   \begin{array}{ll} f_i & = p_i q_{i+1}- p_{i+1} q_i \\
                            f_{i+1} & = p_i q_i + p_{i+1} q_{i+1}
           \end{array} \right\}$ (focus-focus type)

We say that $y$ is a \emph{hyperbolic} singular point if all of its
components are of hyperbolic type, i.e., $f_i = p_iq_i$ for all $i =
1, \hdots,k$ in the above local normal form.

The momentum map ${\bf F}$ gives rise to a singular torus fibration:
by definition, each fiber is a connected component of a level set
(i.e., the preimage of a point in $\bbR^n$) of ${\bf F}$. Regular
fibers of this fibration are Liouville tori, and singular fibers are
those which contain at least one singular point of the system.
Denote by $\cB$ the \emph{base space} of this singular fibration,
with the induced topology from $M$. In general, $\cB$ is a
stratified $n$-dimensional space with an integral affine structure,
and singular points of $\cB$ correspond to the singularities of $\bf
F$ (see \cite{Zung-IntegrableII}). We may consider the Hamiltonian
function $H$ as a function on $\cB$. For each point $z \in \cB$,
denote by $N_z$ the corresponding fiber of the system. If $z$ is a
regular point, i.e. $N_z$ is a Liouville torus, then we will say
that \emph{$H$ satisfies the Kolmogorov condition at $z$} if there
is a local integral affine coordinate system $(I_1,\hdots,I_n)$ near
$z$ on $\cB$ (i.e., a local system of action variables) such that
$\det (\partial^2 H / \partial I_i\partial I_j ) (z) \neq 0$.

Consider now a singular fiber $N = N_x$ of the system. We say that
$N$ is a nondegenerate hyperbolic singularity of corank $k$, if the
following two conditions are satisfied (see \cite{ZungA-L} and
Chapter 9 of \cite{BolsinovFomenko-IntegrableBook}):

1) Each point of $N$ is either regular, or nondegenerate hyperbolic
singular of corank smaller or equal to $k$, and there is at least
one nondegenerate hyperbolic singular point of corank $k$ on $N$.

2) The \emph{non-splitting condition} (which was called the
``topological stability condition" in \cite{ZungA-L}): there is a
neighborhood of $N$ in $M$, such that when we restrict ${\bf F}$ to
this neighborhood, then the set of its singular values in $\bbR^n$
(i.e., the local bifurcation diagram) is a union of $k$ local
transversal (i.e., in generic position) smooth hypersurfaces
intersecting at ${\bf F}(N)$

Consider now such a hyperbolic singularity $N_x$ of corank $k$
(where $x$ denotes the corresponding singular point on the base
space $\cB$). Denote by $y$ a hyperbolic singular point $y$ of
corank $k$ in $N$. By Vey-Eliasson theorem, there is a local
symplectic coordinate system $(p_1,q_1,\hdots, p_n,q_n)$ in which
the $n$ functions $f_1 = p_1, \hdots, f_{n-k} = p_{n-k}, f_{n-k+1} =
p_{n-k+1}q_{n-k+1}, \hdots, f_n = p_nq_n$ are first integrals of the
system. We will make the following assumption about $H$:

3) \emph{H is really nondegenerate hyperbolic} at $N$, in the sense
that when writing $H$ as a function of $n$ variables
$f_1,\hdots,f_n$, we have ${\partial H \over
\partial f_{n-k+i} } (y) \neq 0$ for all $i=1,\hdots,k$. (In other words, the
eigenvalues of the reduced linearized Hamiltonian system of $H$ are
all non-zero real numbers). Remark that this condition does not
depend on the choice of the corank $k$ point $y$ in $N$.

According to the \emph{topological decomposition theorem} for
nondegenerate singularities \cite{ZungA-L}, there is a neighborhood
$({\cU}(N_x), {\cL})$ of $N_x$ in $M$ together with the singular
torus foliation $\cL$ of the system, which is diffeomorphic to an
almost direct product of corank 1 hyperbolic singularities. In other
words, we may write
\begin{equation}
\label{eqn:Decomposition} ({\cU}(N_x), {\cL})  \stackrel{\rm
diffeo}{\simeq} \left({\Bbb T}^{n-k} \times D^{n-k} \times ({\cU}_1,
{\cL}_1) \times ... \times ({\cU}_k, \cL_k) \right) / \Gamma,
\end{equation}
where  ${\Bbb T}^{n-k} \times D^{n-k}$ denotes a trivial fibration
by $(n-k)$-dimensional tori over an $(n-k)$-dimensional disk, each
$({\cU}_i, {\cL}_i)$ is a 2-dimensional surface together with a
singular circle fibration given by the level sets of a Morse
function with one hyperbolic singular level set (there may be many
singular points on the singular level set), $\Gamma$ is a finite
group which acts freely and component-wise on the product (its
action on $D^{n-k}$ is trivial).

Note that the above direct decomposition is \emph{not} symplectic,
i.e. the symplectic form $\omega$ on ${\cU}(N_x)$ cannot be written
as a direct sum of the symplectic forms on the components in
general. However, according to \cite{ZungA-L}, $({\cU}(N_x), {\cL})$
admits a partial system of action-angle variables. In particular,
there is a system of $(n-k)$ action functions $(I_1,\hdots,I_{n-k})$
defined in $({\cU}(N_x), {\cL})$ which gives rise to a locally free
Hamiltonian $\bbT^{n-k}$-action which preserves the system.

The singular point $y$ of corank $k$ in $N$ projects to a singular
point $\hat{y}$ of corank $k$ in $({\cU}_1, {\cL}_1) \times ...
\times ({\cU}_k, \cL_k)$ modulo $\Gamma$. The set $P = ({\Bbb
T}^{n-k} \times D^{n-k} \times \{\hat{y}\}) /\Gamma$ is a symplectic
submanifold in $M$, called the \emph{center manifold} of the system
through $y$. The restriction of our integrable Hamiltonian system to
this center manifold $P$ is a regular integrable Hamiltonian system
with action functions $I_1,\hdots,I_{n-k}$. Our last condition on
$H$ is the following:

4) If $k < n$ then  the restriction $H_P$ of $H$ to the center
manifold $P = ({\Bbb T}^{n-k} \times D^{n-k} \times \{\hat{y}\})
/\Gamma$ satisfies the Kolmogorov condition at the
$(n-k)$-dimensional torus containing $y$ on $P$: $ \det (\partial^2
H_P /
\partial I_i\partial I_j )_{i,j \leq n-k} (y) \neq 0$.

Remark that the above condition does not depend on the choice of the
corank $k$ point $y$ in $N$, and can also be reparaphrased as
follows: $x$ lies on a $(n-k)$-dimensional stratum $S$ in $\cB$ with
a local system of affine coordinates $I_1,\hdots,I_{n-k}$, and we
require that the restiction $H_S$ of $H$ to this stratum $S$ satisfy
the condition $ \det (\partial^2 H_S /
\partial I_i\partial I_j )_{i,j \leq n-k} (x) \neq 0$.

Finally, changing the momentum map in ${\cU}(N_x)$ without changing
the associated singular torus fibration of the system, we can assume
that the following condition on the momentum map is satisfied:

5) $F_1 = I_1,\hdots, F_{n-k} = I_{n-k}$ are action functions, and
for each $i=1,\hdots,k$, $F_{n-k+i}$ is a Morse function on the
component $({\cU}_i, {\cL}_i)$ in the decomposition
(\ref{eqn:Decomposition}) which gives rise to the singular fibration
${\cL}_i$, is equal to zero on the singular fiber of $({\cU}_i,
{\cL}_i)$, and is invariant under the action of $\Gamma$.

\begin{thm}
\label{thm:Kcondition} Consider a smooth integrable Hamiltonian
system with Hamiltonian function $H$ and with a proper momentum map
${\bf F} = (F_1,\hdots,F_n)$ on a $2n$-dimensional symplectic
manifold $(M,\omega)$, which admits a nondegenerate hyperbolic
singularity $N_x$ of corank $k$ ($1 \leq k \leq n$). Assume that $H$
satisfies the conditions 3) and 4) above, i.e., $H$ is really
nondegenerate hyperbolic  at $N_x$, and satisfies the Kolmogorov
condition on a corresponding local $2(n-k)$-dimensional center
manifold. Assume moreover that the momentum map ${\bf F}$ has been
chosen in such a way that it satisfies the above condition 5). Then
we have the following asymptotic formula:
\begin{equation}
\label{eqn:AsymptoticDetBis} \det ( \partial^2 H / \partial
I_i\partial I_j )_{i,j \leq n}(z) = {g(z) \over \prod_{i=1}^k F_{n
-k + i}(z)(\ln F_{n -k + i}(z))^3},
\end{equation}
where $z$ denotes a regular point on the base space $\cB$ of the
system near $x$, $g(z)$ is a smooth first integral in a connected
component of a neighborhood of $x$ in $\cB$ minus the singular part,
such that the limit $\lim_{z \to x} g(z)$ exists and is different
from zero, and $I_1,\hdots, I_n$ is a system of action functions (in
the regular connected component which contains $z$). In particular,
$\det ( \partial^2 H / \partial I_i\partial I_j )_{i,j \leq n}(z)
\neq 0$ for any regular point $z$ which lies in a sufficiently small
neighborhood of $x$ in $\cB$.
\end{thm}

Recall that a particular (and maybe most practical) case of the
above theorem is when the hyperbolic singularity $N_x$ is of corank
$n$, i.e. when it contains a fixed point. In that case, the only
additional condition (Condition 3) on $H$ is that the eigenvalues of
the linear part of the Hamiltonian vector field of $H$ at a fixed
point on $N_x$ are all different from zero.

\section{Asymptotic formula for action functions}

We will keep the notations of the previous section. Consider a
hyperbolic singularity $N_x$ of corank $k$. Recall that in a
neighborhood of $N_x$ there are $n-k$ regular actions functions
$I_1,\hdots,I_{n-k}$. In this section we will write down an
asymptotic formula for the remaining (singular) $k$ action functions
in a complete system of action functions.

Remark that the actions functions change by an affine
transformation, and the determinant $\det (\partial^2 H / \partial
I_i\partial I_j )$ changes by a non-zero multiplicative constant,
when we replace $\cU(N_x)$ by a finite covering of it and lift the
system to that finite covering. So without loss of generality, and
for convenience, from now on we will assume that our singularity
$N_x$ is of direct product type, i.e. the finite group $\Gamma$ in
the decomposition (\ref{eqn:Decomposition}) is trivial:
\begin{equation}
\label{eqn:DirectDecomposition} ({\cU}(N_x), {\cL}) \stackrel{\rm
diff}{\simeq} {\Bbb T}^{n-k} \times D^{n-k} \times ({\cU}_1,
{\cL}_1) \times ... \times ({\cU}_k, \cL_k)
\end{equation}

We will assume that the momentum map has been chosen in such a way
that it satisfies condition 5) of the previous section, i.e. $F_i =
I_i$ for $1 \leq i \leq n-k$ and $F_{n-k+i}$ is a Morse function on
$({\cU}_i, {\cL}_i)$ which gives rise to the fibration ${\cL}_i$ for
$1 \leq i \leq k$, and such that $F_{n-k+i} = 0$ on the singular
fiber of $({\cU}_i, {\cL}_i)$.

Consider a regular point $z$ near $x$ in the base space $\cB$, so
that the Liouville torus $N_z$ lies in the neighborhood ${\cU}(N_x)$
of $N_x$. We can view the momentum map as a map from $\cB$ to
$\bbR^n$. Without loss of generality, we can assume that
$$F_1(x) = \hdots = F_n(x) = 0$$ and $$F_{n-k+1}(z) > 0, \hdots,
F_n(z) > 0.$$ Under the direct decomposition
(\ref{eqn:DirectDecomposition}), we have
\begin{equation}
N_z = \bbT^{n-k}(z) \times S_1(z) \times \hdots \times S_k(z),
\end{equation}
 where $\bbT^{n-k}(z)$ is a fiber in $\bbT^{n-k} \times D^{n-k}$ and
each $ S_i(z)$ is a regular circle fiber in $({\cU}_i, {\cL}_i)$ on
which $F_{n-k+i}$ is constant and positive.

Denote by $C$ the closure of intersection of the local regular
stratum which contains $z$ in $\cB$ with the base of $({\cU}(N_x),
{\cL})$ (i.e. the image of the projection of ${\cU}(N_x)$ to $\cB$).
The set $C$ may be identified with a neighborhood of $0$ of the
``corner'' set $\{ (F_1,\hdots,F_n) \in \bbR^n | F_{n-k+1} \geq 0 \
\forall \ i=1,\hdots,k \}$, with local coordinates
$(F_1,\hdots,F_n)$.

On the interior of $C$ we have two different coordinate systems: the
momentum coordinate system $(F_1,\hdots,F_n)$, and an action
coordinate system $(I_1,\hdots,I_n)$, where $I_1, \hdots,I_{n-k}$
are action variables mentioned above (recall that $F_1 = I_1,
\hdots, F_{n-k} = I_{n-k}$), and each $I_{n-k+i}$ ($i=1,\hdots,k$)
is an action variable defined as follows:
$$ I_{n-k+i} (z) = \int_{\gamma_i(z)} \theta,$$
where $\theta$ is a primitive of the symplectic form $\omega$ in
${\cU}(N_x)$ (i.e., $d\theta = \omega$), and $\gamma_i(z)$ is the
1-cycle on $N_z$ which is represented by $S_i(z)$. On $C$, the
Hamiltonian $H$ is a smooth function of the variables
$F_1,\hdots,F_n$, but $I_{n-k+i}$ ($i=1,\hdots,k$) are not. The
following proposition about the asymptotic behavior of the $k$
singular action functions $I_{n-k+i}$, viewed as functions of $n$
variables $(F_1,\hdots,F_n)$ on $C$ near the origin, will be the
main ingredient in the proof of Theorem \ref{thm:Kcondition}:

\begin{prop} \label{prop:asymptotic}
With the above notations and assumptions, we have, for
$i=1,\hdots,k$,
$$I_{n-k+i} = \psi_i  F_{n-k+i} \ln F_{n-k+i} + \phi_i,$$
on $C$, where $\psi_i = \psi_i(F_1,\hdots,F_n)$ and $\phi_i = \phi_i
(F_1,\hdots,F_n)$ are smooth functions of $n$ variables
$(F_1,\hdots,F_n)$, and $\psi_i(0,\hdots,0) \neq 0$.
\end{prop}

In particular, the action functions $I_{n-k+i}$ admit a continuous
extension on the boundary of $C$ (because $F_{n-k+i} \ln F_{n-k+i}$
tends to $0$ when $F_{n-k+i}$ tends to $0$).

The above proposition is not a new result: it has been known for
some time to people (e.g., Alexey Bolsinov and Vu Ngoc San
\cite{BolsinovSan-SymplecticInv2007}) who work on symplectic
invariants of integrable Hamiltonian systems, and is a direct
consequence of the theorems of Eliasson \cite{Eliasson-NF1990} and
Miranda and myself \cite{EvaZung-Orbit2004} on the local canonical
normal form of an integrable Hamiltonian system near a nondegenerate
singular point or orbit. Let us sketch here its proof:

For simplicity, first consider the case with $n = k = 1$. In this
case, we have just one first integral $F$, one singular action
function $I$, and up to a constant and a sign, $I(z)$  is equal to
the symplectic area of the region $R(z)$ between the singular fiber
$F^{-1}(0)$ and the regular fiber which contains $z$. Near each
singular point $y_i$ ($i = 1, \hdots,m$, where $m$ is the number of
hyperbolic singular points on the singular fiber) we have a local
symplectic coordinate system $(p_i,q_i)$ in which the local
fibration of the system is given by $p_iq_i = constant$. Denote by
$D_i = \{- \epsilon < p_i <\epsilon, - \epsilon < q_i <\epsilon\}$
charts around $y_i$ chosen small enough so that they don't
intersect. The region $R(z)$ can be cut into ``singular pieces''
$R_i(z) = R_z \cap D_i$ and the rest $\hat{R}(z) = R(z) \setminus
\cup_{i} R_i(z)$ (at least one of the singular pieces $R_i(z)$ is
non-empty). The symplectic area of $\hat{R}(z)$ is a smooth function
with respect to $F$, while the symplectic area of each non-empty
singular piece $R_i(z)$ is of the type $\psi F \ln F + \phi$ where
$\psi$ and $\phi$ are smooth with respect to $F$, with $\psi(0) <
0$. Summing up these symplectic area gives us the desired formula
for $I(z)$.

The general (higher dimensional and higher corank) case is the same.
The main idea is to cut a loop on $N_z$ which represents the 1-cycle
$\gamma_{i}$ into several pieces; the integral of the primitive form
$\theta$ over those pieces which pass nearby singular points will
contribute singular terms of the type $\psi_i F_{n-k+i} \ln
F_{n-k+i}$.

\section{Proof of Theorem \ref{thm:Kcondition}}
\label{section:proof}

We will work under the assumptions of Theorem \ref{thm:Kcondition},
and with the notations introduced in the previous sections. Denote
by $\Gamma = (\Gamma_1,\hdots,\Gamma_n)$ the frequency map, where
$\Gamma_i ={\partial H / \partial I_i}$. We will first view
$(\Gamma_i)$ as a map of $n$ variables $(F_i)$ and find an
asymptotic formula for $\det (\partial \Gamma_i / \partial F_j)$,
and then deduct from that asymptotic formula the desired asymptotic
formula for $ \det (\partial^2 H / \partial I_i\partial I_j )$.

To simplify the formulas, we will use the following notations: by
$(smooth)$ we  mean a function on $C$ which is smooth with respect
to the variables $(F_1,\hdots,F_n)$ (they must be smooth also on the
boundary of $C$), by $(smooth*)$ we mean a smooth function which
moreover does not vanish at the origin, by $(small)$ a continuous
function of the variables $(F_1,\hdots,F_n)$ which vanishes at the
origin, by $(smoothsmall)$ a function which is both $(smooth)$ and
$(small)$, by $(continuous)$ a continuous function on $C$, and by
$(continuous*)$ a continuous function which does not vanish at the
origin.

It follows from Proposition \ref{prop:asymptotic} that we have:
\begin{equation}
{\partial I_{n-k+i} \over \partial F_{n-k+i}} = (smooth*). \ln
F_{n-k+i} + (smooth)
\end{equation}
(for $i\leq k$), and
\begin{equation}
{\partial I_{n-k+i} \over \partial F_{j}} = (smooth). F_{n-k+i} \ln
F_{n-k+i} + (smooth)
\end{equation}
(for $j \neq n-k+i;\  i \leq k; \ j \leq n$.)

Since $I_i = F_i \ \forall i \leq n-k$, we obviously have
\begin{equation}
{\partial I_{i} \over \partial F_{i}} = 1
\end{equation}
and
\begin{equation}
{\partial I_{i} \over \partial F_{j}} = 0
\end{equation}
for all $i \leq n-k, j \neq i, j \leq n$.

The asymptotic behavior (near the origin) of the entries of the
matrix $\left( {\partial I_i \over \partial F_j}
\right)^{i=1,\hdots, n}_{j=1,\hdots, n}$ are given by the above
formulas. Let us now write down the asymptotic formulas for the
entries of the inverse matrix $\left( {\partial F_i \over \partial
I_j} \right)^{i=1,\hdots, n}_{j=1,\hdots, n}$. Direct computations
show that:

\begin{equation}
\det \left( {\partial I_i \over \partial F_j} \right)^{i=1,\hdots,
n}_{j=1,\hdots, n} = (smooth*).\prod_{j=1}^k \ln F_{n-k+j} +
(l.o.t.)
\end{equation}
where $(l.o.t.)$ (lower order terms) means terms of the following
types:
\\ $(smoothsmall).\prod_{j=1}^k \ln F_{n-k+j}$ for some $i$, and
$(smooth).\prod_{j \in \Delta} \ln F_{n-k+j}$ where $\Delta$ is a
proper subset of $\{1,\hdots,k\}$;

\begin{equation}
{\partial F_{n-k+i} \over \partial I_{n-k+i}}  = {(smooth*).\prod_{j
\neq i} \ln F_{n-k+j} + (l.o.t.) \over (smooth*).\prod_{j} \ln
F_{n-k+j} + (l.o.t.)}
\end{equation}
(for $i \leq k$), where $l.o.t.$ in $(smooth*).\prod_{j \neq i} \ln
F_{n-k+j} + (l.o.t.)$ mean terms of the following types:
$(smoothsmall).\prod_{j\neq i} \ln F_{n-k+j}$, and
$(smooth).\prod_{j \in \Delta} \ln F_{n-k+j}$ where $\Delta$ is a
proper subset of $\{1,\hdots,k\} \setminus \{i\}$ (i.e. terms of
smaller order than $\prod_{j \neq i} \ln F_{n-k+j}$);

\begin{equation}
{\partial F_{n-k+i} \over \partial I_{n-k+s}}  =
{(smooth).F_{n-k+i}\prod_{j \neq s} \ln F_{n-k+j} +
(smooth).\prod_{j \neq s,i} \ln F_{n-k+j} + (l.o.t.) \over
(smooth*).\prod_{j} \ln F_{n-k+j} + (l.o.t.)}
\end{equation}
(for $s \neq i$), and

\begin{equation}
{\partial F_{n-k+i} \over \partial I_{t}}  =
{(smooth).F_{n-k+i}\prod_{j} \ln F_{n-k+j} + (smooth).\prod_{j \neq
i} \ln F_{n-k+j} + (l.o.t.) \over (smooth*).\prod_{j} \ln F_{n-k+j}
+ (s.o.t)}
\end{equation}
(for $t \leq n-k$). The reader may have noticed that our (partial)
ordering of the terms is generated by
$$ \ln F_{n+k-i} \succ (smooth*) \succ (smoothsmall).$$

The above formulas together with the formula $\Gamma_i =
\sum_{j=1}^n {\partial H \over \partial F_j}.{\partial F_j \over
\partial I_i}$ (for $i=1\hdots,n$) give rise to:

\begin{equation}
\Gamma_t  = {\partial H \over \partial F_t} + {\sum_i (smooth).
F_{n-k+i}\prod_{j} \ln F_{n-k+j} + \sum_{i}(smooth).\prod_{j \neq i}
\ln F_{n-k+j} + (l.o.t.) \over (smooth*).\prod_{j} \ln F_{n-k+j} +
(l.o.t)}
\end{equation}
(for $t \leq n-k$) and

\begin{equation}
\Gamma_{n-k+i}  = {(smooth*).\prod_{j \neq i} \ln F_{n-k+j} +
(l.o.t.) \over (smooth*).\prod_{j} \ln F_{n-k+j} + (l.o.t.)}
\end{equation}
(for $i \leq k$).

The above asymptotic formulas for the frequency map lead directly to
the following formulas:

\begin{equation}
{\partial \Gamma_t \over \partial F_s} = {\partial^2H \over \partial
F_t \partial F_s} + (small)
\end{equation}
(for $t,s \leq n-k$);

\begin{equation}
{\partial \Gamma_{n-k+i} \over \partial F_s} =  (small)
\end{equation}
(for $i \leq k, s \leq n-k$);

\begin{equation}
{\partial \Gamma_{n-k+j} \over \partial F_{n-k+j}} = { (continuous*)
\over F_{n-k+j} (\ln  F_{n-k+j})^2 }
\end{equation}
(for $j \leq k$);

\begin{equation}
{\partial \Gamma_{n-k+i} \over \partial F_{n-k+j}} = { (continuous)
\over F_{n-k+j} (\ln  F_{n-k+j})^2 (\ln F_{n-k+i})} =  { (small)
\over F_{n-k+j} (\ln  F_{n-k+j})^2 }
\end{equation}
(for $i,j \leq k, i \neq j$); and

\begin{equation}
{\partial \Gamma_t \over \partial F_{n-k+j}} = { (continuous) \over
F_{n-k+j} (\ln  F_{n-k+j})^2 }
\end{equation}
(for $j \leq k, t \leq n-k$).

In turn, the above asymptotic formulas for the entries of the matrix
$\left( {\partial \Gamma_i \over \partial F_j} \right)^{i \leq n}_{j
\leq n}$ imply that the leading term in the asymptotic expansion of
the determinant $\det \left( {\partial \Gamma_i \over \partial F_j}
\right)$ is of the form
\begin{equation}
 \det \left({\partial^2H (x) \over \partial F_s \partial
F_t}\right)_{s,t \leq n-k}. \prod_{j \leq k} { (continuous*)\over
F_{n-k+j} (\ln F_{n-k+j})^2},
\end{equation}
 where $ \det \left({\partial^2H (x)
\over \partial F_s \partial F_t}\right)_{s,t \leq n-k} \neq 0$ by
our hypothesis, so we can write
\begin{equation}
\det \left( {\partial \Gamma_i \over \partial F_j} \right) =
{(continuous*) \over \prod_{j \leq k} { F_{n-k+j} (\ln
F_{n-k+j})^2}}
\end{equation}

It follows from the asymptotic formula for the matrix $\left(
{\partial I_i / \partial F_j} \right)$ shown earlier in this section
that we have
\begin{equation}
\det \left( {\partial I_i \over \partial F_j} \right) =
{(continuous*)  \prod_{j \leq k} \ln F_{n-k+j}}
\end{equation}

The last two formulas, together with the fact that $$ \det
(\partial^2 H / \partial I_i\partial I_j ) = \det (\partial \Gamma_i
/ \partial I_j) = \det ( \partial \Gamma_i /
\partial F_s ) / \det \left( \partial I_j /
\partial F_s \right) $$
give us the asymptotic formula
\begin{equation}
\det ( \partial^2 H / \partial I_i\partial I_j ) = {(continuous*)
\over \prod_{i=1}^k F_{n -k + i} (\ln F_{n -k + i})^3}
\end{equation}
on $C$. The theorem is proved.


\begin{thebibliography}{50}
\parskip0.0cm

\small{

\baselineskip0.35cm

\bibitem{Arnold}
V. I. Arnold, Mathematical methods of classical mechanics,
Springer-Verlag 1978.

\bibitem{BolsinovFomenko-IntegrableBook}
A. V. Bolsinov and A. T. Fomenko, Integrable Hamiltonian systems:
Geometry, topology, classification, 2004, CRC.

\bibitem{BolsinovSan-SymplecticInv2007}
Alexey Bolsinov and Vu-Ngoc San, {\it Symplectic equivalence for
integrable systems with common action integrals}, in preparation
(2007).

\bibitem{DullinSan-Twist2004}
Holger Dullin and Vu-Ngoc San, {\it Vanishing twist near focus-focus
points} Nonlinearity 17 (2004), 1777-1785.

\bibitem{Eliasson-NF1990}
Hakan Eliasson, {\it Normal forms for {H}amiltonian systems with
{P}oisson commuting
  integrals---elliptic case}, Comment. Math. Helv. \textbf{65} (1990), no.~1,
  4--35; and the PhD thesis, 1984.

\bibitem{Horozov}
Emil Horozov, {\it Perturbations of the spherical pendulum and
Abelian integrals}, J. reine angew. Math., 408 (1990), 114-135.

\bibitem{Knorrer}
Horst Kn\"orrer, {\it Singular fibres of the momentum mapping for
integrable Hamiltonian systems}, J. Reine Angew. Math., 355 (1985),
67-107.

\bibitem{Kolmogorov}
A. N. Kolmogorov, Selected works, Vol. 1 (V.M. Tikhomirov ed.),
Cluwer Acad. Publ., 1991.

\bibitem{EvaZung-Orbit2004}
Eva Miranda and Nguyen Tien Zung, {\it Equivariant normal form for
nondegenerate singular orbits of integrable Hamiltonian systems},
Ann. Sci. Ecole Norm. Sup. 37 (2004), No. 6, 819-839.

\bibitem{Moser}
Jurgen Moser, Stable and Random Motions in Dynamical Systems, Ann.
Math. Studies 77, Princeton 1973.

\bibitem{Mineur-AA1935}
Henri Mineur, \emph{{Sur les systemes mecaniques admettant $n$
integrales
  premieres uniformes et l'extension a ces systemes de la methode de
  quantification de Sommerfeld}}, C. R. Acad. Sci., Paris \textbf{200} (1935),
  1571--1573 (French).

\bibitem{Mineur-AA1937}
Henri Mineur, \emph{{Sur les systemes mecaniques dans lesquels
figurent des
  parametres fonctions du temps. Etude des systemes admettant $n$ integrales
  premieres uniformes en involution. Extension a ces systemes des conditions de
  quantification de Bohr-Sommerfeld.}}, Journal de l'Ecole Polytechnique, Série
  III, 143ème année (1937), 173--191 and 237--270.

\bibitem{Rink-Focus2004}
Bob Rink, {\it A Cantor set of tori with monodromy near a
focus-focus singularity}, Nonlinearity 17 (2004), Number 1, 347-356.

\bibitem{Russmann-KAM2001}
Helmut Rüssmann, {\it Invariant tori in non-degenerate nearly
integrable Hamiltonian systems}, Regular and Chaotic Dynamics 6
(2001), 119-204.

\bibitem{Vey-Separable1978}
Jacques Vey, {\it Sur certaines syst\`emes dynamiques s\'eparables},
Amer. J. Math., 100 (1978), 591-614.

\bibitem{ZungA-L} Nguyen Tien Zung,
{\it Symplectic topology of integrable Hamiltonian systems, I:
Arnold-Liouville with singularities}, Compositio Math., 101 (1996),
179-215.


\bibitem{Zung-IntegrableII} Nguyen Tien Zung,
{\it Symplectic topology of integrable Hamiltonian systems, II:
Topological classification}, Compositio Math., 138 (2003), No. 2,
125-156.

}

\end{thebibliography}
\end{document}